\documentclass[11pt]{article}
\newtheorem{conjecture}{Conjecture}

\newtheorem{corollary}[conjecture]{Corollary}

\newtheorem{lemma}[conjecture]{Lemma}

\newtheorem{proposition}[conjecture]{Proposition}

\newtheorem{remark}[conjecture]{Remark}
\newtheorem{theorem}[conjecture]{Theorem}

\newenvironment{proof}{\noindent {\bf Proof} \hspace{.1cm}}{\hfill ${\bf QED}$ \\ \vspace{.15cm}}

\title{Uniformly exponential growth and mapping class groups of surfaces}

\author{James W. Anderson, Javier Aramayona and Kenneth J. Shackleton\footnote{Supported by an EPSRC studentship}}
\date{17 August, 2005}

\begin{document}

\maketitle

\begin{abstract}
\noindent
We show that the mapping class group of an orientable finite type surface has uniformly exponential growth, as well as various closely related groups.  This provides further evidence that mapping class groups may be linear.



\medskip
\noindent
MSC 20F65 (primary), 20F38 (secondary)


\end{abstract}

The purpose of this short note is to demonstrate that the mapping class group of an orientable finite type surface has uniformly exponential growth.  This result is new for surfaces of genus at least one, with the exception of the closed surface of genus two (see below).  This should be viewed as further evidence that such mapping class groups may be linear, as many linear groups are known to have uniformly exponential growth by recent work of Eskin, Mozes, and Oh \cite{emo}.  Alternatively, this can be viewed as removing a possible avenue for showing that such mapping class groups are not linear.  In this sense, our work is similar in spirit to the recent work of Brendle and Hamidi-Tehrani \cite{brendle-hh}. We go on to show that closely related groups of homotopy classes of homeomorphisms of surfaces, as well as analogous groups of automorphisms of free groups, also have uniformly exponential growth. We remark that, while the linearity of most surface mapping class groups is an open question, most automorphism groups of free groups are known not to be linear. Specifically, as noted in Bendle and Hamidi-Tehrani  \cite{brendle-hh}, it is known that ${\rm Aut}(F_n)$ is not linear for $n\ge 3$ and ${\rm Out}(F_n)$ is not linear for $n\ge 4$, whereas both ${\rm Aut}(F_2)$ and ${\rm Out}(F_2)$ are linear.


We begin by reviewing some basic definitions.  For a survey of exponential growth and uniformly exponential growth, we refer to the article of de la Harpe \cite{de la harpe} and the references contained therein.


Let $G$ be a finitely generated group, and let $S$ be a finite generating set for $G$.  The {\em length} $\ell_S(g)$ of an element $g\in G$ is the minimum integer $k$ so that $g$ can be expressed as $g = s_{i_1}\cdots s_{i_k}$, where each $s_{i_j}\in S\cup S^{-1}$.  We define the length of the identity element of $G$ to be $0$.  Let \[ B_S(n) =\{ g\in G\: |\: \ell_S(g)\le n\} \] be the (closed) ball of radius $n$ about the identity element in $G$, and let $| B_S(n)|$ be the number of elements of $G$ in $B_S(n)$.


The {\em exponential growth rate} $\omega(G,S)$ of $G$ with respect to $S$ is defined to be
\[ \omega(G,S) =\lim_{n\rightarrow\infty} \sqrt[n]{|B_S(n)|}. \] (This limit exists, due to the submultiplicativity $\ell_S(gh)\le \ell_S(g) \ell_S(h)$ of the length function on $G$.)
The group $G$ has {\em exponential growth} if $\omega(G,S) >1$ for some (and hence for every) finite generating set $S$. Note that if $G$ has a free subgroup of rank $2$, then $G$ has exponential growth, though not conversely.

We can also remove the dependence on particular generating sets by considering \[ \omega(G) =\inf_S \omega(G,S), \] where the infimum is taken over all finite generating sets $S$ of $G$.  The group $G$ is said to have {\em uniformly exponential growth} if $\omega(G) >1$.   For examples of groups of exponential growth that do not have uniformly exponential growth, see Wilson \cite{wilson}.


We make use of the following Proposition from de la Harpe \cite{de la harpe}.

\begin{proposition} [Proposition $2.3$ of de la Harpe \cite{de la harpe}] If $G$ is a finitely generated group and if $G'$ is a quotient of $G$, then $\omega(G)\ge \omega(G')$.
\label{lower-bound}
\end{proposition}


We use the following result of Shalen and Wagreich \cite{shalen-wagreich} to show the uniformly exponential growth for certain groups related to the mapping class group (see Corollary \ref{punctures} and Corollary \ref{boundary}).

\begin{lemma} [Corollary $3.6$ of Shalen and Wagreich \cite{shalen-wagreich}] Let $G$ be a finitely generated group and let $H$ be a finite index subgroup of $G$.  If $H$ has uniformly exponential growth, then $G$ has uniformly exponential growth.
\label{commensurability}
\end{lemma}

\begin{remark} To the best of our knowledge, it is not yet known whether the converse of Lemma \ref{commensurability} holds. It seems the main difficulty lies in singling out an extended generating set for $G$ from one for $H$ to give uniform embeddings on Cayley graphs. In this note, such issues represent only a minor inconvenience. A positive answer would prove uniformly exponential growth to be a commensurability invariant.  We conjecture that this should be true.

\end{remark}


The main tool we use is the following result, which should be viewed as an immediate corollary of the recent work of Eskin, Mozes, and Oh \cite{emo}, in which they demonstrate that finitely generated subgroups of ${\rm GL}(n,{\bf C})$ have uniformly exponential growth if and only if they have exponential growth.

\begin{theorem} Let $G$ be a finitely generated group, and suppose that for some $n\ge 2$
there exists a homomorphism $\rho: G\rightarrow {\rm GL}_n({\bf C})$ whose image $\rho(G)$ is not virtually nilpotent.  Then, $G$ has uniformly exponential growth. 
\label{emo-corollary}
\end{theorem}

\begin{proof} Since $\rho(G)$ is finitely generated and not virtually nilpotent, Corollary $1.4$ of Eskin, Mozes, and Oh \cite{emo} yields $\omega(\rho(G)) > 1$ (that is, $\rho(G)$  has uniformly exponential growth). Since $\rho(G)$ is the homomorphic image of $G$, it is a quotient of $G$. Proposition \ref{lower-bound} yields that $\omega(G)\ge \omega(\rho(G)) > 1$, and so $G$ has uniformly exponential growth.
\end{proof}

Note that it is known that virtually nilpotent groups have polynomial growth (and conversely, by work of Gromov \cite{gromov}), whereas groups containing a free group of rank at least two have exponential growth.  


Let $\Sigma$ be a closed orientable surface of genus $g\ge 1$, and let $P$ be a finite set of marked points on $\Sigma$, where $n = |P|\ge 0$.  The {\em mapping class group} ${\cal M}(\Sigma,P) ={\cal M}_{g,n}$ is the set of homotopy classes of orientation preserving homeomorphisms $f: \Sigma\rightarrow \Sigma$ for which $f(P) =P$, where the homotopies are required to keep each point of $P$ fixed. Note that the elements of ${\cal M}_{g,n}$ can permute the points of $P$. We refer to the survey article of Ivanov \cite{ivanov} as our main reference for the mapping class group and its properties.

There is a natural surjective homomorphism from ${\cal M}_{g,n}$ to ${\rm Symm}(n)$, the symmetric group on $n$ letters, given by restricting the action of ${\cal M}_{g,n}$ to $P$. The kernel ${\cal PM}_{g,n}$ of this homomorphism is  the {\em pure mapping class group}, which is the subgroup of ${\cal M}_{g,n}$ that fixes every element of $P$. Note that ${\cal PM}_{g,n}$ is a subgroup of finite index in ${\cal M}_{g,n}$.


Now suppose $n=|P|>0$ and let $p\in P$. By forgetting the marked point $p$ we obtain that every homeomorphism $f:\Sigma \rightarrow \Sigma$ fixing $P$ pointwise induces a homeomorphism $f':\Sigma\rightarrow \Sigma$ fixing $P\setminus \{p\}$ pointwise. In this way we get a surjective homomorphism ${\cal PM}_{g,n} \rightarrow {\cal PM}_{g,n-1}$; in particular, we see that ${\cal PM}_{g,n}$ homomorphically surjects onto the mapping class group $\mathcal{PM}_{g,0}={\cal M}_{g,0} = {\cal M}_g$ of the surface $\Sigma$ with no marked points.

The {\em extended mapping class group} $\mathcal{M}^{\pm}_{g,n}$ is the group of homotopy classes of all homeomorphisms of $\Sigma$ fixing $P$ setwise, and is a degree $2$ extension of $\mathcal{M}_{g,n}$.

We can make similar definitions in the case that $\Sigma$ has non-empty boundary.  Let $\Sigma$ be a compact orientable surface of genus $g\ge 1$ with non-empty boundary, and let $P$ be a finite set of $n =|P| \ge 0$ marked points in the interior of $\Sigma$.  Suppose that $\partial \Sigma$ has $m\ge 1$ components.  The {\em mapping class group} ${\cal M}_{g,n,m}$ is the group of homotopy classes of orientation preserving homeomorphisms of $\Sigma$.  The {\em pure mapping class group} ${\cal PM}_{g,n,m}$ is the subgroup of ${\cal M}_{g,n,m}$ consisting of those elements that permute neither the components of $\partial \Sigma$ nor the elements of $P$. As before, ${\cal PM}_{g,n,m}$ is a subgroup of finite index in ${\cal M}_{g,n,m}$. The {\em extended mapping class group} ${\cal M}^{\pm}_{g,n,m}$ is the group of homotopy classes of all homeomorphisms of $\Sigma$, and is a degree $2$ extension of ${\cal M}_{g,n,m}$.

We note that for a compact orientable surface $\Sigma$ with non-empty boundary, there is a natural surjective homomorphism from $\mathcal{PM}_{g,n,m}$ to $\mathcal{PM}_{g,n,0} =\mathcal{PM}_{g,n}$, obtained by gluing discs to all the boundary components of $\Sigma$ and extending the homeomorphisms of $\Sigma$ across these discs; this is discussed in detail in Theorem 2.8.C of Ivanov \cite{ivanov}.

We are now ready for the main result of this note.

\begin{theorem} For $g\ge 1$, the groups ${\cal M}_{g}$ and ${\cal
M}^{\pm}_{g}$ have uniformly exponential growth.
\label{main-mcg}
\end{theorem}


\begin{proof} Let $\Sigma$ be a closed surface of genus $g\ge 1$ and consider the
mapping class group ${\cal M}_g$ of $\Sigma$. We show the result
for ${\cal M}_{g}$; Lemma \ref{commensurability} then yields
the result for ${\cal M}^{\pm}_{g}$, which is a degree 2 extension
of ${\cal M}_{g}$ as we noted before. The action of ${\cal M}_g$
on $\pi_1(\Sigma)$ descends to a linear action of ${\cal M}_g$ on
$H_1(\Sigma,  {\bf Z}) =\pi_1(\Sigma)/[\pi_1(\Sigma),\pi_1(\Sigma)]$, yielding a surjective homomorphism $\rho: {\cal
M}_g\rightarrow {\rm Sp}(2g, {\bf Z})$. Since ${\rm Sp}(2g, {\bf
Z})$ contains ${\rm Sp}(2, {\bf Z}) \cong {\rm SL}(2, {\bf Z})$ as
a subgroup and since ${\rm SL}(2,{\bf Z})$ contains a ${\bf Z}\ast
{\bf Z}$ subgroup, we see that $\rho({\cal M}_g) = {\rm Sp}(2g,
{\bf Z})$ is not virtually nilpotent.  Hence, by Theorem
\ref{emo-corollary}, we see that ${\cal M}_{g}$ has uniformly
exponential growth.

[To see that  ${\rm Sp}(2g, {\bf Z})$ contains ${\rm Sp}(2, {\bf
Z}) \cong {\rm SL}(2, {\bf Z})$, recall that ${\rm Sp}(2g, {\bf
Z})$ is defined to be the group of $2g\times 2g$-matrices
preserving a non-degenerate, skew-symmetric bilinear form.  If we
take the form ${\bf q}({\bf x}, {\bf y}) = x_1 y_2 - x_2 y_1 +
\cdots + x_{2g-1} y_{2g} - x_{2g}
 y_{2g-1}$, then any block diagonal matrix of the form
\[ \left( \begin{array}{cc} A & {\bf 0}_{2, 2g-2} \\ {\bf 0}_{2g-2, 2} & I_{2g-2} \end{array} \right) \]
preserves ${\bf q}$ and hence lies in
${\rm Sp}(2g, {\bf Z})$, where $A$ lies in the group ${\rm Sp}(2, {\bf
Z})$ preserving the quadratic form ${\bf q}_0({\bf x}, {\bf y})
=x_1 y_2 - x_2 y_1$, ${\bf 0}_{s,t}$ is the $s\times t$ matrix of
zeroes, and $I_{k}$ is the $k\times k$ identity matrix.] 

\end{proof}

We saw before that there is a surjective homomorphism from ${\cal PM}_{g,n}$ to ${\cal PM}_{g}= {\cal M}_{g}$, and that ${\cal PM}_{g,n}$ has finite index in ${\cal M}_{g,n}$.  Proposition \ref{lower-bound} and Lemma \ref{commensurability} then yield the following corollary.

\begin{corollary} For $g \geq 1$ and $n>0$, the groups ${\cal
PM}_{g,n}$, ${\cal M}_{g,n}$ and ${\cal M}^{\pm}_{g,n}$ have uniformly exponential growth.
\label{punctures}
\end{corollary}

We also saw that, in the case of surfaces with boundary, there is
a homomorphism ${\cal PM}_{g,n,m}$ to ${\cal PM}_{g,n}$, and that
${\cal PM}_{g,n,m}$ has finite index in ${\cal M}_{g,n,m}$. We thus get the following result, again using Proposition \ref{lower-bound} and Lemma
\ref{commensurability} .

\begin{corollary} For $g \geq 1$ and $m,n>0$, the groups ${\cal
PM}_{g,n,m}$, ${\cal M}_{g,n,m}$ and ${\cal M}^{\pm}_{g,n,m}$ have uniformly exponential growth.
\label{boundary}
\end{corollary}

For $g=0$, the methods we use here do not apply.  We note that when $n\ge 4$, it is a remarkable result of Bigelow \cite{bigelow} and Krammer \cite{krammer} that each ${\cal M}_{0,n}$ is linear and hence has uniformly exponential growth.  For $g=2$ and $n=0$, it is a result of Bigelow and Budney \cite{bigelow-budney} and of Korkmaz \cite{korkmaz} that ${\cal M}_{2,0}$ is linear and hence has uniformly exponential growth.




We note in closing that analogous results hold for the {\em automorphism group} ${\rm Aut}(F_n)$ of the free group $F_n$ of rank $n$ and for the {\em outer automorphism group} ${\rm Out}(F_n)$, the quotient of ${\rm Aut}(F_n)$ by the group of inner automorphisms.

\begin{theorem} For $n\ge 2$, the groups ${\rm Aut}(F_n)$ and ${\rm Out}(F_n)$ have uniformly exponential growth.
\label{main-out}
\end{theorem}

\begin{proof} We begin with ${\rm Out}(F_n)$.  The action of ${\rm Out}(F_n)$ on $F_n$ descends to a linear action of ${\rm Out}(F_n)$ on $F_n/[F_n, F_n] ={\bf Z}^n$, yielding a
surjective homomorphism $\rho: {\rm Out}(F_n)\rightarrow {\rm GL}(n,{\bf Z})$.  Since ${\rm GL}(n,{\bf Z})$ is finitely generated and is not virtually nilpotent (as ${\rm SL}(2,{\bf
Z})\subset {\rm SL}(n,{\bf Z})\subset {\rm GL}(n,{\bf Z})$), Theorem \ref{emo-corollary} yields that ${\rm Out}(F_n)$ has uniformly exponential growth.

That ${\rm Aut}(F_n)$ has uniformly exponential growth follows from Proposition \ref{lower-bound}, since ${\rm Out}(F_n)$ is a quotient of ${\rm Aut}(F_n)$.
\end{proof}



\noindent 
{\bf{Acknowledgements.}} Parts of this paper were written during visits of the second author to the Bernoulli Center, Lausanne, and the third author to the Tokyo Institute of Technology. Both authors wish to express their gratitude to these institutions.  The first author would also like to thank Gareth Jones for helpful conversations.

{\footnotesize J. W. Anderson\\School of Mathematics\\University of Southampton\\Southampton SO17 1BJ\\England\\j.w.anderson@maths.soton.ac.uk}

{\footnotesize J. Aramayona (corresponding author)\\Mathematics Institute\\University of Warwick\\Coventry CV4 7AL\\England\\jaram@maths.warwick.ac.uk}

{\footnotesize K. J. Shackleton\\School of Mathematics\\University of Southampton\\Southampton SO17 1BJ\\England\\k.j.shackleton@maths.soton.ac.uk}

\end{document}